%% file: main.tex
\documentclass[11pt]{article}

\usepackage[a4paper,margin=0.8in]{geometry}
\input{header}

\usepackage{soul}
\usepackage{color}

\title{Open Problem: Two Riddles in Heavy-Ball Dynamics}

\author{
    Baptiste Goujaud
    \thanks{INRIA Saclay, Palaiseau, France. \texttt{baptiste.goujaud@inria.fr}
    \newline
    \setlength{\parindent}{14pt}
    \indent
    \setlength{\parindent}{15pt}
    Work done while at CMAP, École Polytechnique, Institut Polytechnique de Paris, France.
    }
    \and 
    Adrien Taylor\thanks{INRIA, D.I. École Normale Supérieure, CNRS, PSL Research University, Paris, France. \texttt{adrien.taylor@inria.fr}} 
    \and 
    Aymeric Dieuleveut\thanks{CMAP, École Polytechnique, Institut Polytechnique de Paris, France. \texttt{aymeric.dieuleveut@polytechnique.edu}}
}

\begin{document}

\maketitle
\vspace{-1em}
\begin{abstract}
This short paper presents two open problems on the widely used Polyak's Heavy-Ball algorithm. The first problem is the method's ability to exactly \textit{accelerate} in dimension one exactly. The second question regards the behavior of the method for parameters for which it seems that neither a Lyapunov nor a cycle exists.

For both problems, we provide a detailed description of the problem and elements of an answer.
\end{abstract}

\subsection*{Preamble: What constitutes an Open Problem?} An open problem in mathematics is typically characterized as an unsolved, well-defined question, on a topic that has attracted interest from the community, and which resolution can bring novel insights to the field. 

The first problem we describe pertains to a straightforward and natural question regarding the long-standing issue of the behavior of the Heavy-Ball method beyond quadratics, a topic on which several contradictory results have been published. While we acknowledge that the impact of the result may be limited beyond its original mathematical curiosity, our complementary goal is to clarify the current state of knowledge on this matter.

Our second problem is more original. It concerns the type of potential dynamics in first-order optimization. Again, the problem admits a simple and self-contained formulation. The community has shown interest in understanding how systematic Lyapunov approaches can be developed.
The behaviors we hope to uncover may be applicable far beyond the heavy-ball case and may thus help us understand first-order methods in more depth. 

\vspace{-1em}
\section{Introduction, background, and open problems}

First-order optimization methods have recently attracted a lot of attention due to their generally low cost per iteration and their practical success in many applications. They are particularly relevant in applications not requiring very accurate solutions, such as in machine learning (see, e.g.,~\cite{bottou2007tradeoffs}). 
In particular, the Heavy-ball (HB) method, proposed by~\citet{polyak1964some}, is a fundamental algorithm in convex optimization and is widely used due to its complexity improvement over simpler existing algorithms such as Gradient Descent ($\GD$). The novelty of the $\HB$ method over $\GD$ is the addition of a momentum term. For convex differentiable function $f: \R^d \to \R$,  for a step size $\gamma$ and a momentum parameter $\beta$, the update writes:
\begin{equation} 
    x_{t+1} = x_t - \gamma \nabla f(x_t) + \beta (x_t - x_{t-1}).
    \label{eq:hb} \tag{$\HB_{\gb}$}
\end{equation}
This momentum enables to obtain an \textit{acceleration} over the set of $L$-smooth $\mu$-strongly-convex quadratic functions (which we denote by $\Qml$). Precisely, while a well-tuned $\GD$ needs $O(\tfrac{L}{\mu}\log(1/\varepsilon))$ iterations to reach an $\varepsilon$ accuracy on every function of $\Qml$, a well-tuned $\HB$ only needs $O(\sqrt{\tfrac{L}{\mu}}\log(1/\varepsilon))$ iterations to achieve the same precision. This improvement in the dependency of the \emph{conditioning} ($L / \mu$) of the objective function is often referred to as the \emph{acceleration} of the momentum-based method.

The behavior of $\HB$ is (asymptotically) optimal on quadratic functions: indeed an optimally tuned HB method corresponds to the asymptotic version of Chebyshev acceleration~\citep[see e.g.,][]{polyak1964some,fischer2011polynomial}. 
On the other hand, the behavior of $\HB$  on the larger class $\Fml$ of $L$-smooth $\mu$-strongly-convex (non-necessarily quadratic) remains incompletely understood.  \citet{ghadimi2015global} proved convergence results, without acceleration. 
Recently,~\citet{goujaud2023provable} closed the long-standing question acceleration on $\Fml$,  by showing that \textbf{no tuning of $\HB$ allows it to reach such an acceleration} in dimension 2 or more.

\begin{figure}[t]
    \begin{subfigure}{.49\linewidth}
        \centering
        \includegraphics[width=\linewidth]{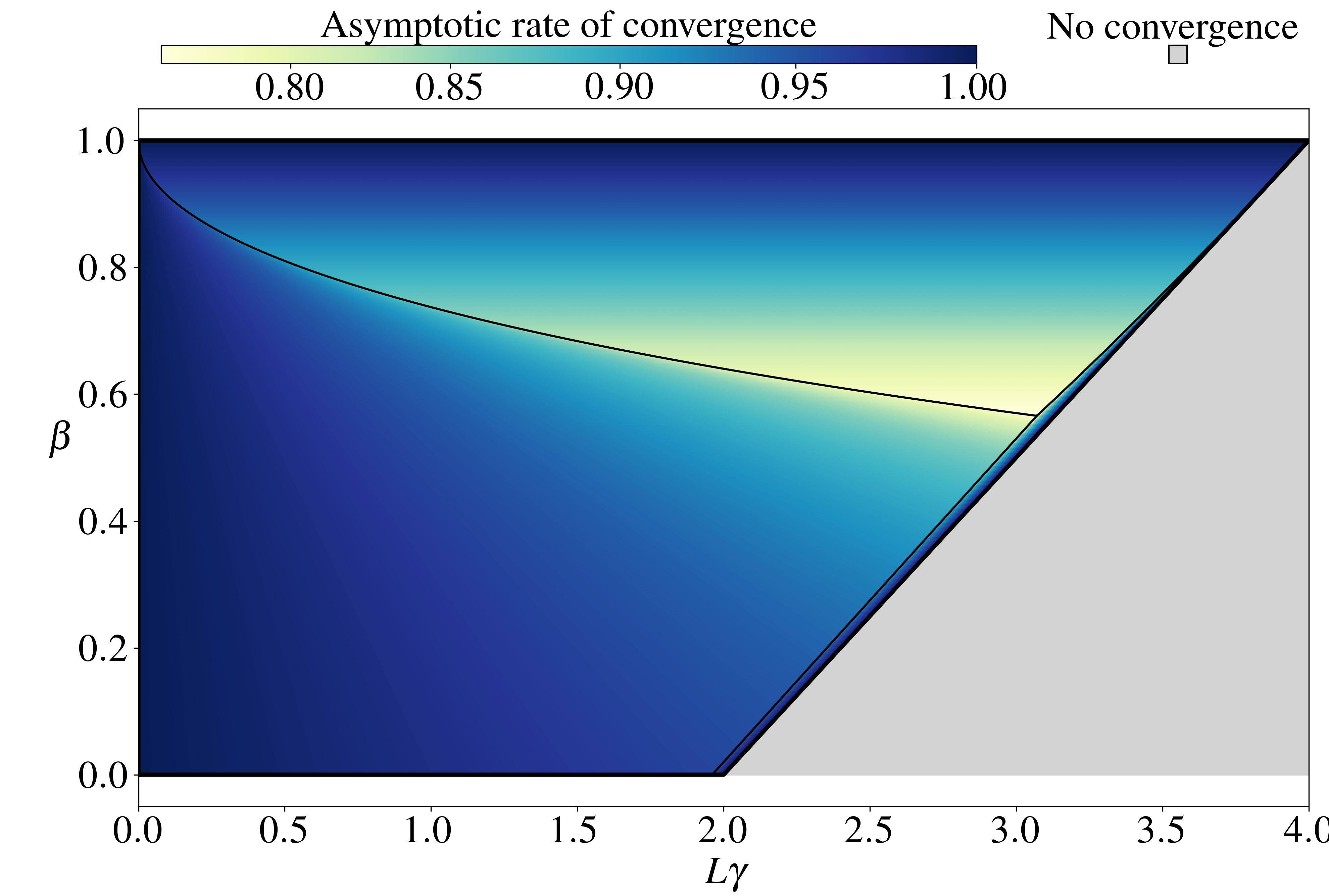}
        \caption{Asymptotic convergence rate on $\Qml$ of HB tuned with the parameters $(\gamma,\beta)$. \label{fig:1a}}
    \end{subfigure}
    \hfill
    \begin{subfigure}{.49\linewidth}
        \centering
        \includegraphics[width=\linewidth]{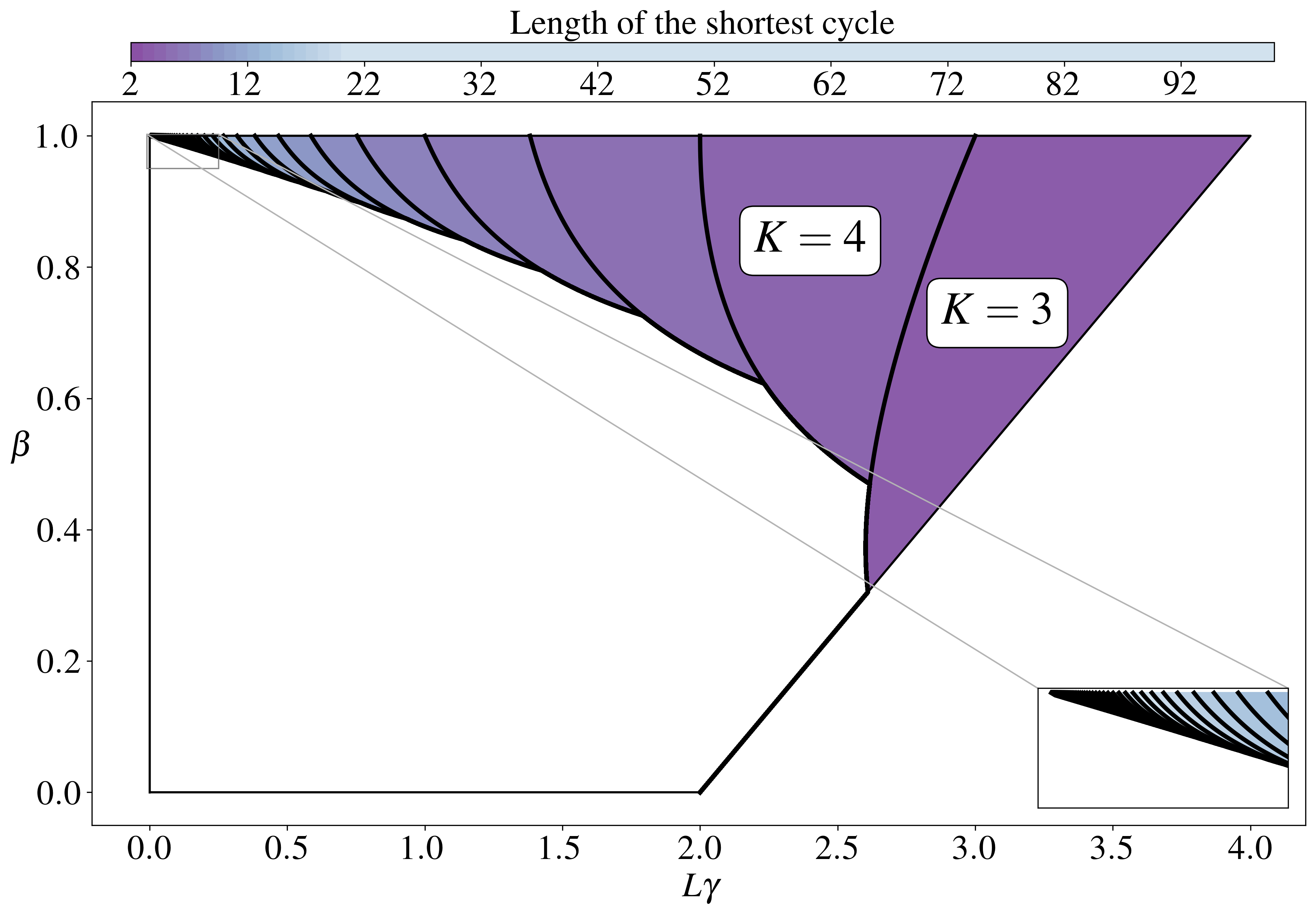}
        \caption{Tunings $( \gamma, \beta)$ of HB, for which a cycle of length $K$ exists on a function of $\Fml$, for $K$ in $\{3, \dots, 25\}$. \label{fig:1b}}
    \end{subfigure}
    \caption{Comparison of HB behaviors on $\Qml$ and $\Fml$.}
\end{figure}

\paragraph{{Sketch of the proof}.} In order to prove such a result, \citet{goujaud2023provable} construct specific smooth strongly-convex functions on which $\HB$ \textit{cycles}, i.e.~there exists an initialization for which the HB trajectory constitutes a cycle of length $K$, and thus never converges. Specifically, for each tuning $(\gamma,\beta)$ of $\HB$ of the purple regions of~\Cref{fig:1b}, such a counter-example is given, in dimension 2, showing that $\HB$ does not converge in the worst-case for those tunings.

\Cref{fig:1a} on the other hand, shows the performance of $\HB$ on the class quadratics functions $\Qml$, which admits a closed-form formula. All the tunings allowing for acceleration (bright region in~\Cref{fig:1a}) on $\Qml$ fall into the purple regions of~\Cref{fig:1b}. That is for any tuning, $\HB$ either converges with a non-accelerated rate on $\Qml$ (and thus on $\Fml$) or does not converge at all on $\Fml$ (see Cor.~3.7 in \citep{goujaud2023provable}).
This analysis leads to two simple remaining open problems.

\begin{center}
    \fbox{\parbox{.9\textwidth}{
        \begin{center}
          \textbf{Open problem 1.}  Does heavy-ball accelerate over \textit{scalar} functions of $\Fml$?
        \end{center}
    }}
\end{center}
Indeed, the counter-examples in~\citep[Theorem 3.5, eq.3]{goujaud2023provable} are functions of 2-dimensional variables ($d=2$). This dimension is specific, as the given counter-examples cycle over the $K$-th \textit{roots of unity}, a shape that arises from problem symmetries. While this proves that  $\HB$ does not accelerate over the set of $d$-dimensional smooth strongly-convex functions for any dimension $d \geq 2$, the behavior of $\HB$ in dimension 1 therefore remains unknown.

\begin{wrapfigure}[14]{r}{0.45\textwidth}
    \begin{center}
        \vspace{-.9cm}
        \includegraphics[width=0.45\textwidth]{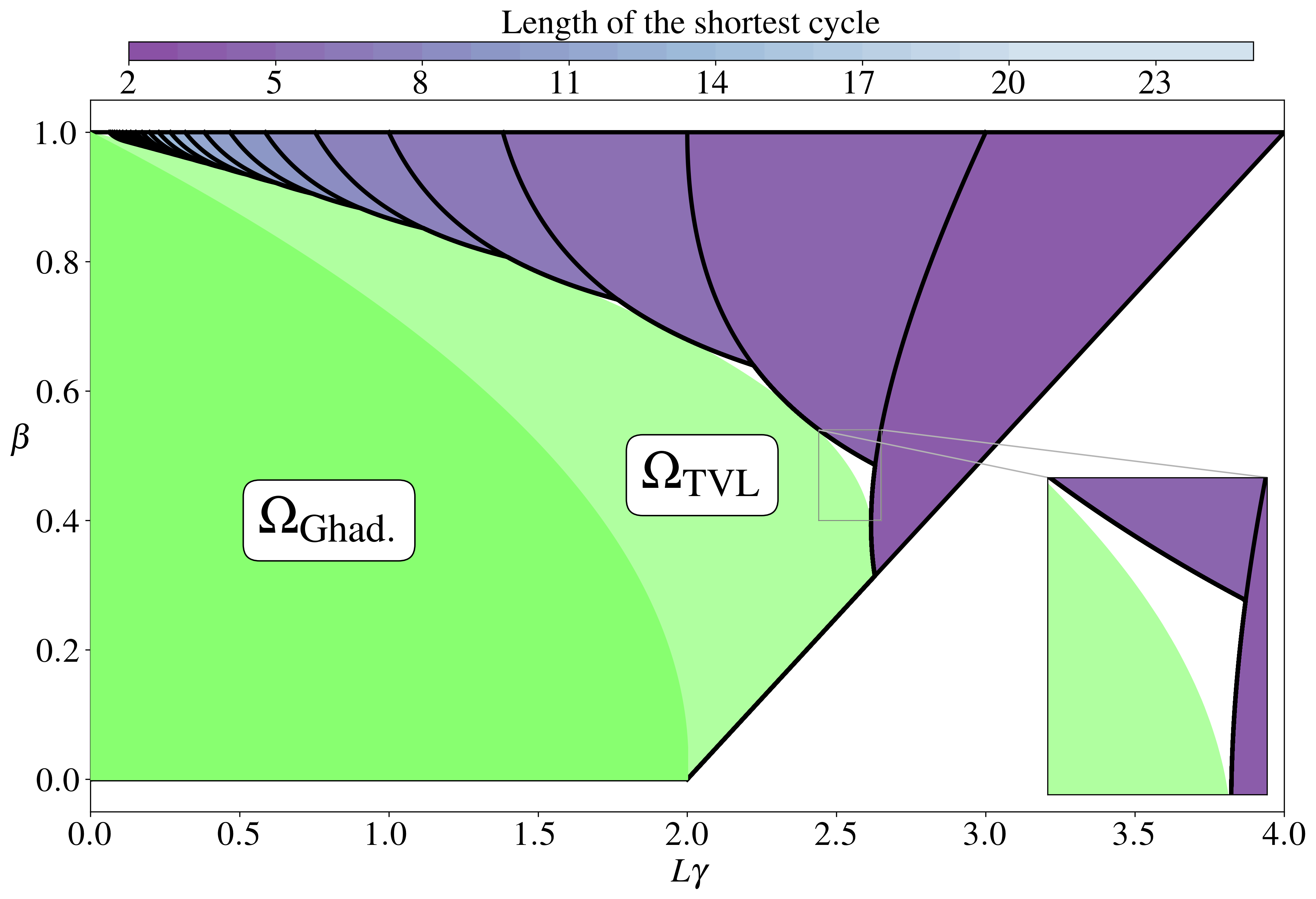}
    \end{center}
    \vspace{-.9cm}
    \caption{$\HB$ behavior as a function of its tuning: Green $\leftrightarrow$ Lyapunov; Purple $\leftrightarrow$ Cycle; Unfilled $\leftrightarrow$ Open Problem 2. \label{fig:2}}
\end{wrapfigure}

\begin{center}
    \fbox{\parbox{.9\textwidth}{
        \begin{center}
            \textbf{Open problem 2.} Do there exist parameters for which HB provably neither admits a Lyapunov function, nor a cycle? 
            If so, what happens then?
        \end{center}
    }}
\end{center}
A surprising fact is that there is strong empirical evidence that the answer to the first part of the question is positive, and the second part is completely open. In~\Cref{fig:2}, the purple region corresponds to the one in which $\HB$ provably cycles on a function of $\Fml$. Moreover, the existence of a cycle can be solved as a linear program: through a \textit{numerical} resolution,  exactly the same set of tunings is obtained.

On the other hand, the green region is where Lyapunov functions are found, ensuring a global convergence of $\HB$ for any tuning in this region. Precisely, the darker green region corresponds to the one found by~\citet{ghadimi2015global} (analytically), while the lighter is obtained numerically, following the procedure described in~\citep{taylor2018lyapunov} using the PEP framework introduced in~\citep{drori2014performance,taylor2017smooth} and the IQC approach described in~\citep{lessard2016analysis}.
Overall, one observes that there exist \textit{unfilled} areas, i.e., tunings for which $\HB$ has neither a Lyapunov function nor a cycle on $\Fml$, and we cannot conclude whether $\HB$ converges or diverges when tuned in those regions. An example of such an area is highlighted in the zoom square in the lower right corner of \Cref{fig:2}.

In the following, we provide elements of answers to these two problems.

\vspace{-1em}
\section{Open problem 1: $\HB$ behavior in dimension 1}
\begin{figure}[h!]
    \begin{subfigure}{.45\linewidth}
        \centering
        \includegraphics[width=\linewidth]{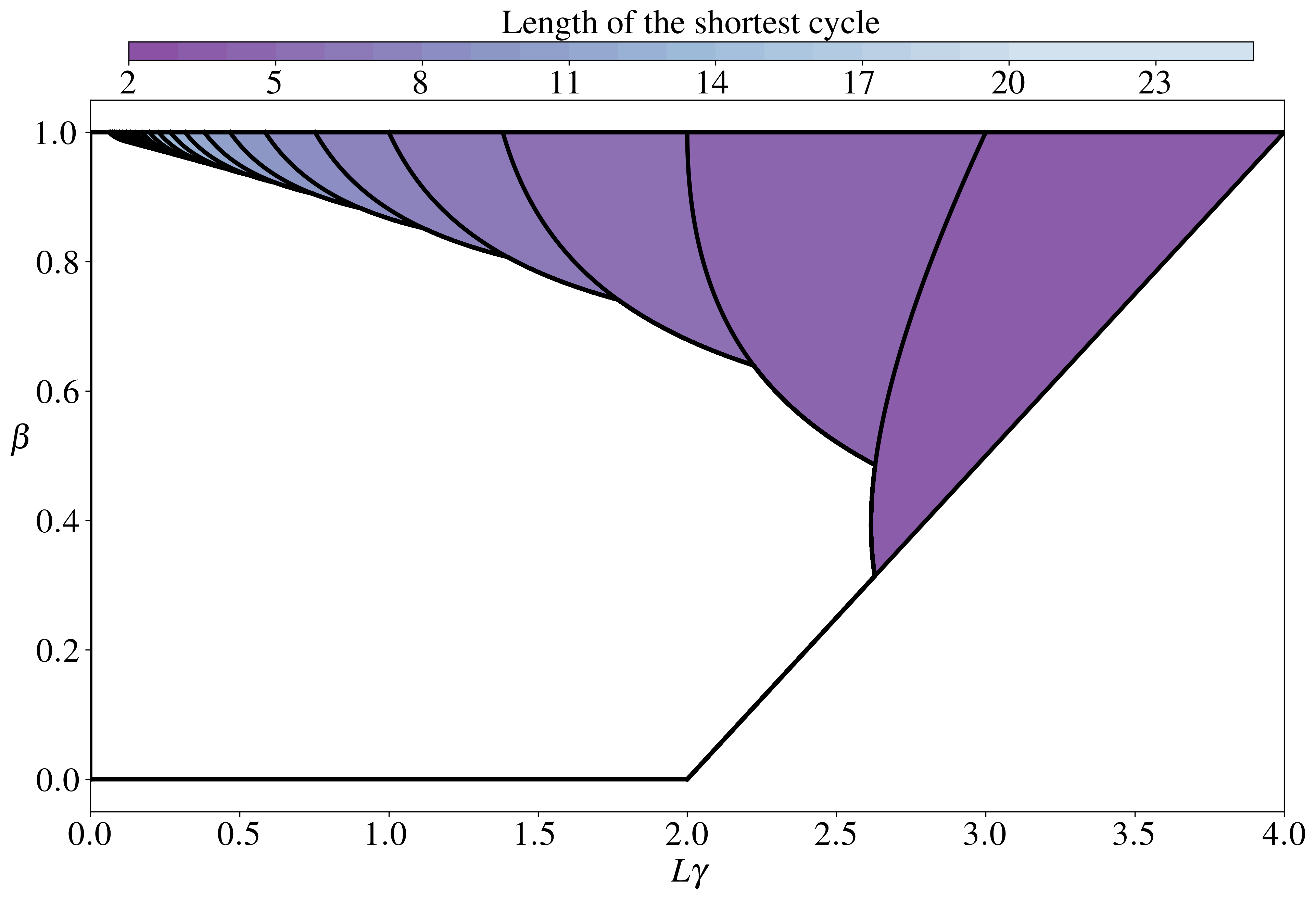}
        \caption{$d>1$ \label{fig:cycle_d_2_inf}}
    \end{subfigure}
    \hfill
    \begin{subfigure}{.45\linewidth}
        \centering
        \includegraphics[width=\linewidth]{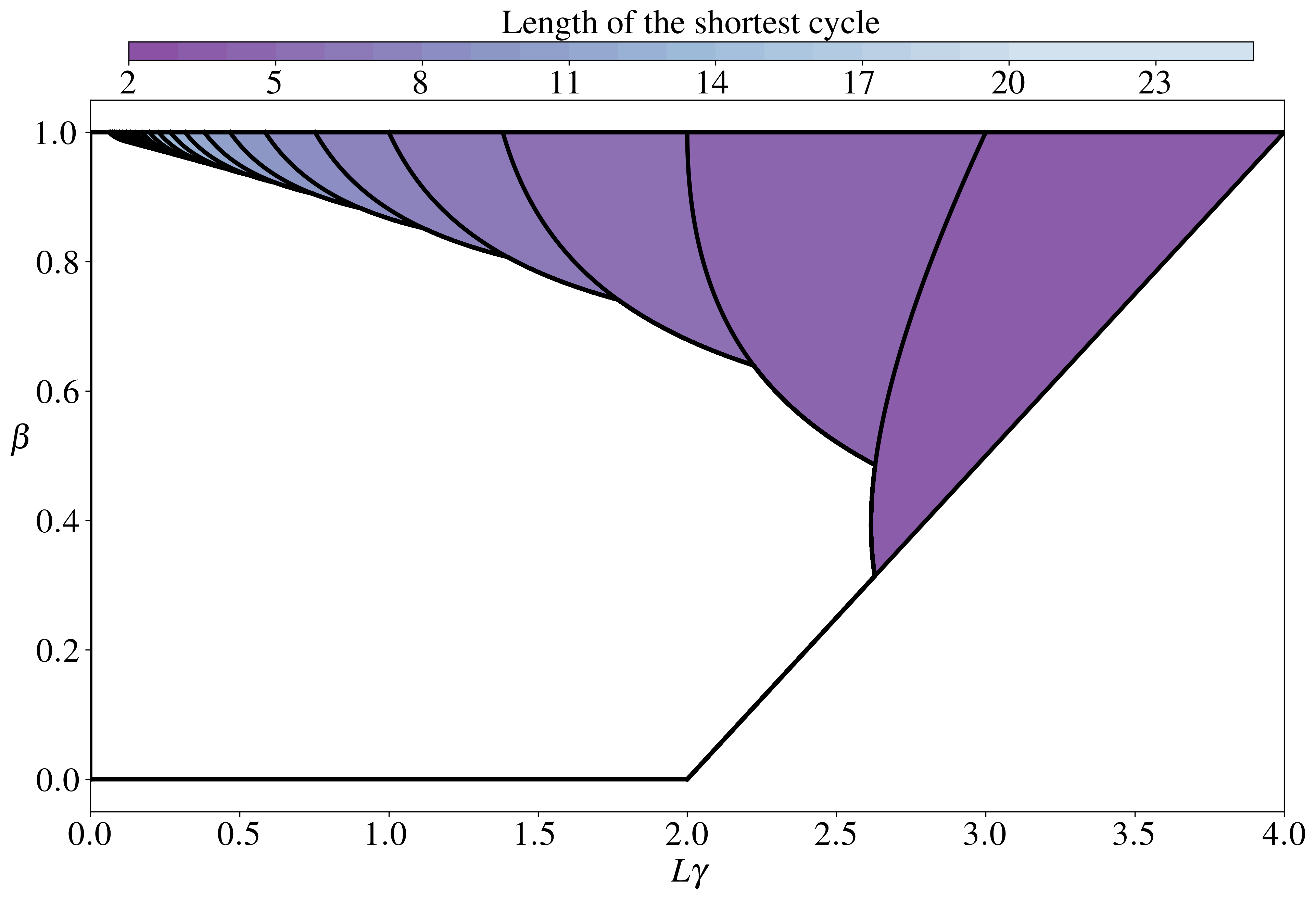}
        \caption{$d=1$ \label{fig:cycle_d_2_1}}
    \end{subfigure}
    \caption{Comparison of cycle regions in different dimensions. Border of $\Orouml$ in black lines. $\Ocyml$ in shades of purples. \textit{Both regions are numerically identical.} \label{fig:dim1}}
    \vspace{-1em}
\end{figure}

In this section, we provide elements that lead us to conjecture that HB also does not accelerate in dimension 1, and elements towards proving that conjecture. In short, we obtain through a numerical approach --detailed hereafter-- the set of tunings (containing at least the one for which \citet{lessard2016analysis} already exhibited a cycle) for which a cycle appears: this set is \textit{indistinguishable} from the set of tunings for which a cycle was established in dimension 2 (that were sufficient to disprove acceleration). This is illustrated in~\Cref{fig:dim1}. While this numerical observation does not constitute a mathematical proof, this leads us to believe that HB cannot accelerate, even in dimension~1.

\paragraph{Numerical approach for finding cycles in dimension 1.}

Let $K \geq 3$. Let $X = (x_0, \cdots, x_{K-1})^\top \in \mathbb{R}^K$ a cycle. A straightforward computation shows that $\HB_{\gb}$ cycles over those iterates when running on a function with gradients $G = (g_0, \cdots, g_{K-1}) \in \mathbb{R}^K$ at points $(x_0, \cdots, x_{K-1}) $, if and only if $G$ and $X$ verify $$G = \frac{[(1+\beta)I - J - \beta J^{-1}]}{\gamma}X,$$ where $J$ is the circular permutation operator, i.e. for all $i$, $J e_i = e_{i+1 \ \mathrm{ mod } \ d}$, for $(e_1, \dots, e_d)$ the canonical eigenbasis.

In dimension 1, verifying that the underlying function belongs to $\Fml$ is equivalent to verifying $\mu \leq \frac{g^{(i+1)} - g^{(i)}}{x^{(i+1)} - x^{(i)}} \leq L$ where the exponents \textit{sort} $X$, i.e. there exists a permutation $\sigma$ such that $X = \sigma (x^{(0)}, \cdots, x^{(K-1)})^\top$ and we define $(g^{(0)}, \cdots, g^{(K-1)})^\top$ with $G = \sigma (g^{(0)}, \cdots, g^{(K-1)})^\top$.
Note $\sigma$ is defined as the permutation that sorts $X$, and it also sorts $G$ if and only if $f$ is convex. In other words, with $\tilde X=(x^{(0)}, \cdots, x^{(K-1)})$ and $\tilde G=(g^{(0)}, \cdots, g^{(K-1)})$, and $\Sigma\in  \mathbb S_{K}$ the permutation matrix associated to $\sigma\in   \mathcal S_{K}$,  there is a cycle of  $\HB_{\gb}$ over $\Fml$ if and only if
\def\arraystretch{1.32}
\begin{equation*}
\begin{array}{ll}
    & \exists X\in \mathbb R^K \ \ \text{s.t. } \quad  (\exists f \in \Fml \quad  \text{s.t. } \HB_{\gb}(f) \text{ cycles on } X  ) \\
    \Leftrightarrow &\exists X \in \mathbb R^K,\ \  \exists \Sigma \in \mathbb  S_{K} \quad \text{s.t. }   \left\{\begin{array}{l}
    \forall i \in [K-1] \quad   0\le \mu \langle \tilde  X, e_{i+1} - e_i\rangle  \leq \langle \tilde G, e_{i+1} - e_i\rangle \leq L \langle\tilde  X, e_{i+1} - e_i\rangle,  \\
    \text{ with } \tilde X = \Sigma^{-1} X  ,   \tilde G = \Sigma^{-1} G, \text{ and  }  G:= \frac{[(1+\beta)I - J - \beta J^{-1}]}{\gamma}X, \\
    \end{array} \right. \\
    \Leftrightarrow   &  \exists \Sigma \in \mathbb  S_{K} \quad \text{s.t. }   \eqref{eq:lpsigma}. \\
\end{array}
\end{equation*}
With, for any  $\Sigma \in \mathbb  S_{K} $,
\begin{equation}
    \label{eq:lpsigma}
    \exists X\in \mathbb R^K,\ \text{s.t. }   \left\{\begin{array}{l}
    \forall i \in [K-1] \quad   0\le  \mu \langle \tilde  X, e_{i+1} - e_i\rangle  \leq \langle \tilde G, e_{i+1} - e_i\rangle \leq L \langle\tilde  X, e_{i+1} - e_i\rangle,  \\
    \text{ with } \tilde X = \Sigma^{-1} X,  \tilde G = \Sigma^{-1} G, \text{ and  }  G:= \frac{[(1+\beta)I - J - \beta J^{-1}]}{\gamma}X. \\ \end{array} \right. \tag{$LP_\Sigma$}
\end{equation}

That is, \textit{for any given permutation} $\sigma$ (or equivalently matrix $\Sigma$), the sub-problem \eqref{eq:lpsigma}  is written as a Linear Program (LP), that can thus be solved efficiently. We use this approach to obtain in~\Cref{fig:cycle_d_2_1} all parameters $\gb$, for which there exists a 1-dimensional cycle for $\HB_{\gb}$ on $\Fml$, for all cycle length up to 6.

\begin{figure}
    \centering
    \begin{subfigure}[t]{.5\linewidth}
        \centering
        \includegraphics[width=.6\linewidth]{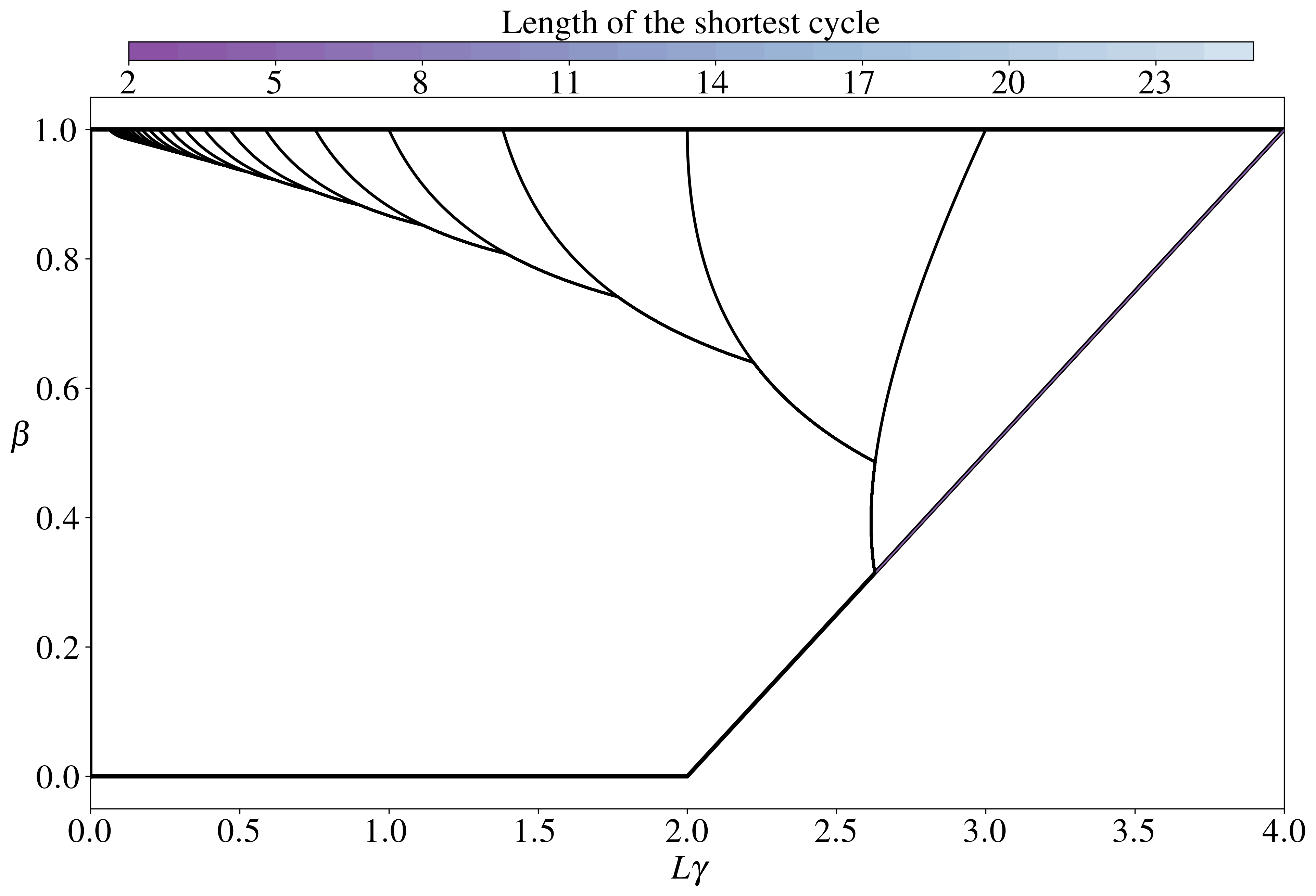}
        \caption{%
        \(\displaystyle 
          \sigma \;=\; 
          \begin{pmatrix}
            0 & 1 & 2 & 3 \\
            0 & 1 & 2 & 3
          \end{pmatrix}
          , \,
          \sigma \;=\; 
          \begin{pmatrix}
            0 & 1 & 2 & 3 \\
            0 & 2 & 1 & 3
          \end{pmatrix}
        \)
        \label{fig:K4_1234}}
    \end{subfigure}\hfill
    \begin{subfigure}[t]{.5\linewidth}
        \centering
        \includegraphics[width=.6\linewidth]{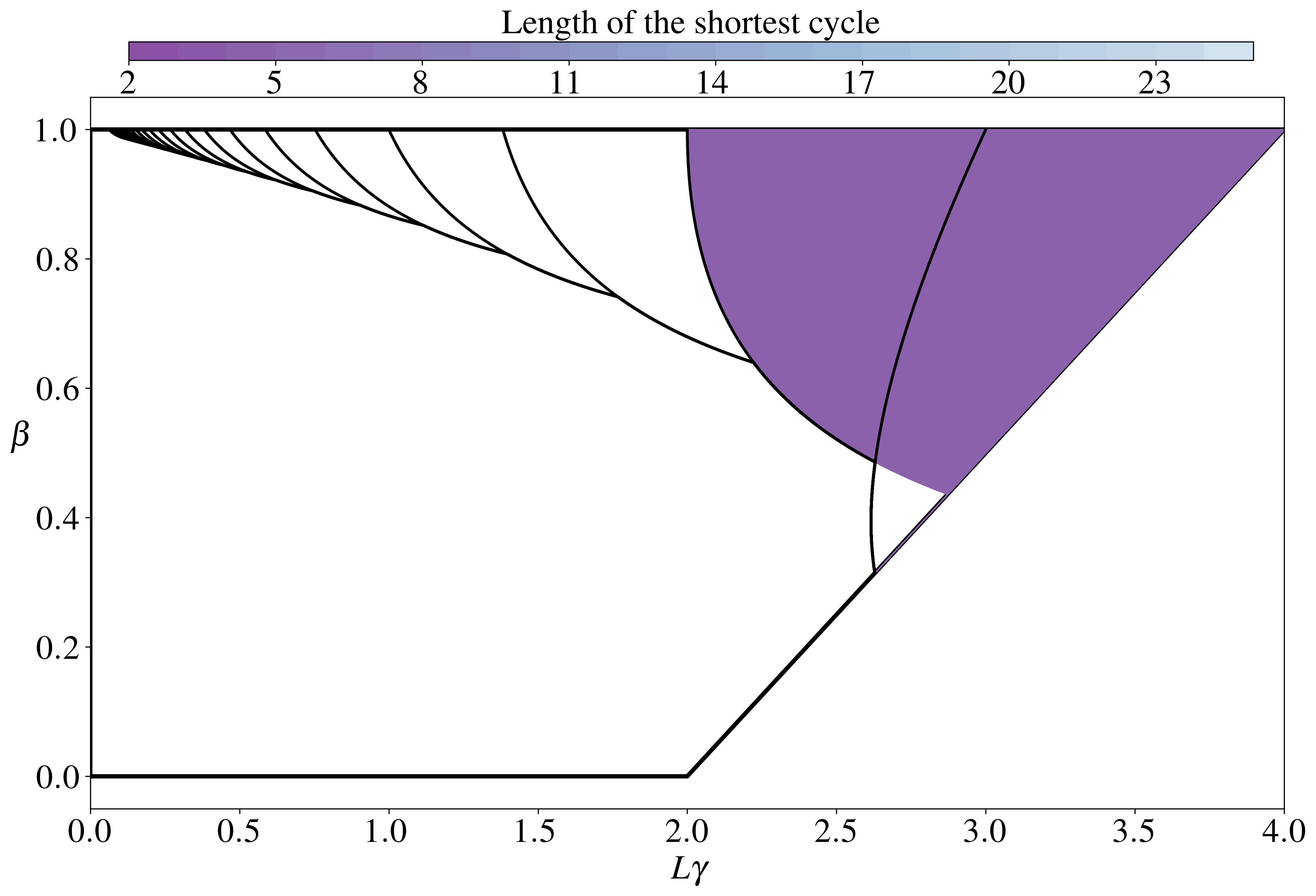}
        \caption{%
        \(\displaystyle 
          \sigma \;=\; 
          \begin{pmatrix}
            0 & 1 & 2 & 3 \\
            0 & 1 & 3 & 2
          \end{pmatrix}
          , \,
          \sigma \;=\; 
          \begin{pmatrix}
            0 & 1 & 2 & 3 \\
            0 & 2 & 3 & 1
          \end{pmatrix}
        \)
        \label{fig:K4_1243}}
    \end{subfigure}

    \caption{For $K=4$, $(\mu,L)$ fixed, the set of points $\gb$ such that~\eqref{eq:lpsigma} admits a feasible point, depending on $\sigma$. \label{fig:cycles_dep_on_sig_K=4}}
\end{figure}

\begin{figure}
    \centering
    \begin{subfigure}[t]{.3\linewidth}
        \centering
        \includegraphics[width=\linewidth]{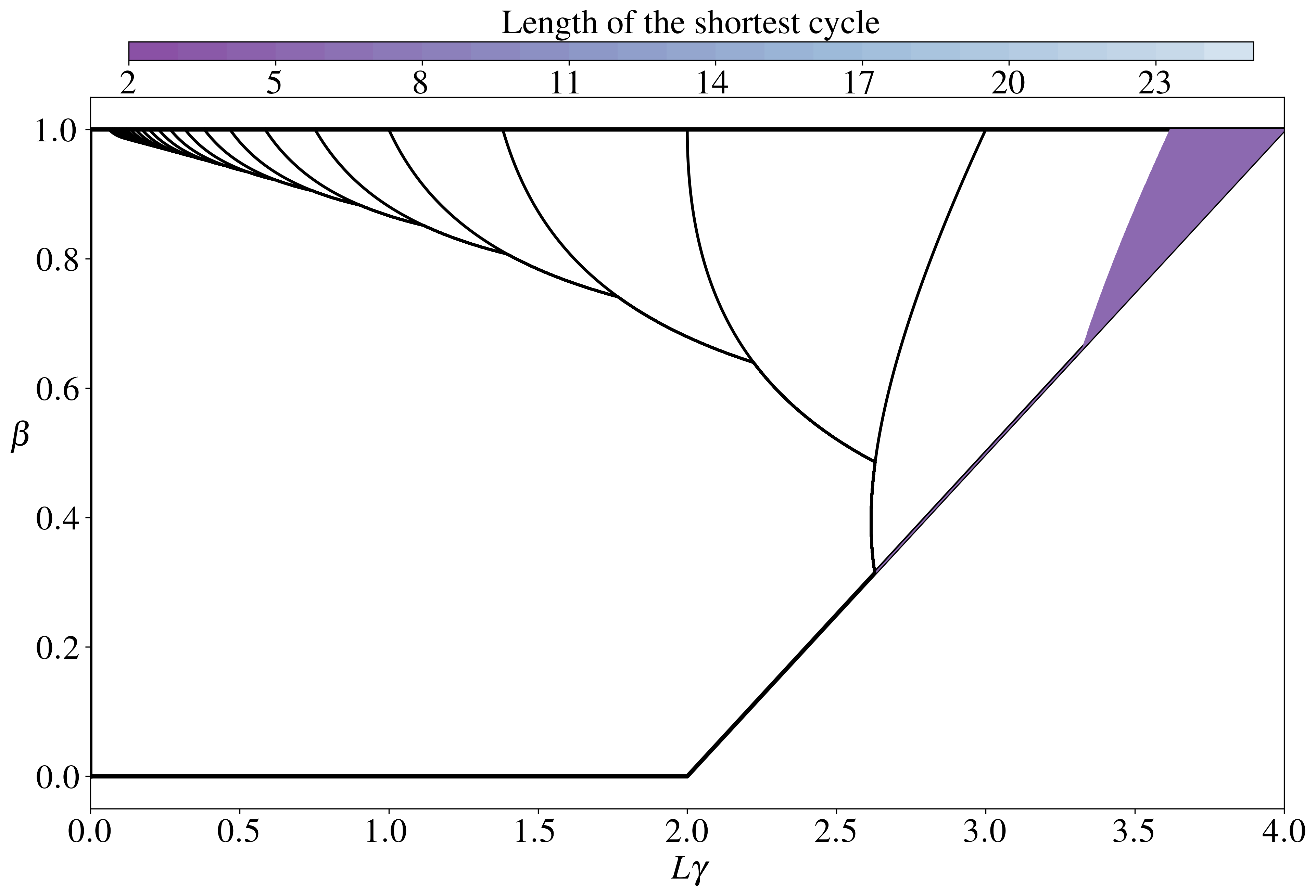}
        \caption{%
        \(\displaystyle 
          \sigma \;=\; 
          \begin{pmatrix}
            0 & 1 & 2 & 3 & 4 \\
            0 & 3 & 2 & 1 & 4
          \end{pmatrix}
        \)
        \label{fig:K5_14325}}
    \end{subfigure}\hfill
    \begin{subfigure}[t]{.3\linewidth}
        \centering
        \includegraphics[width=\linewidth]{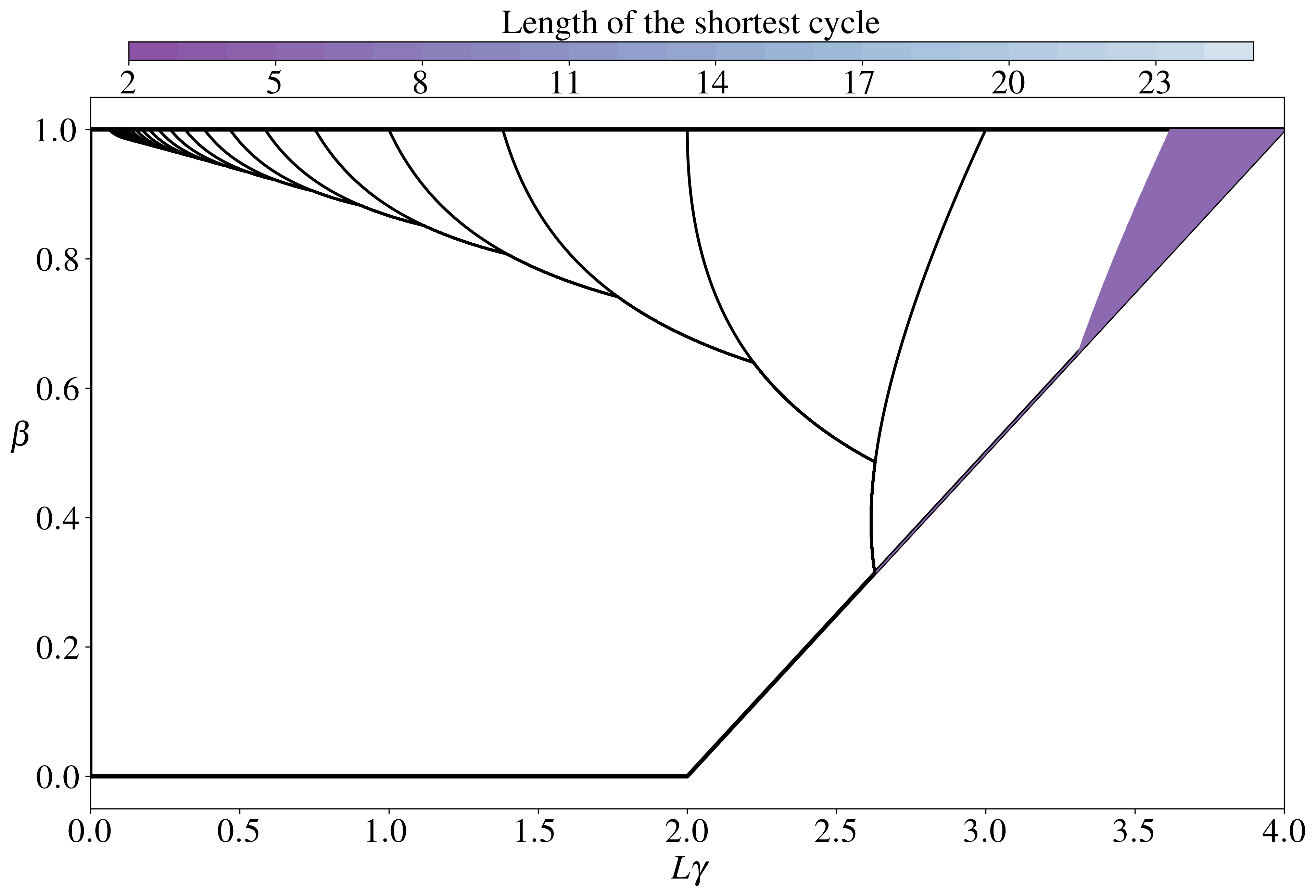}
        \caption{%
        \(\displaystyle 
          \sigma \;=\; 
          \begin{pmatrix}
            0 & 1 & 2 & 3 & 4 \\
            0 & 4 & 1 & 3 & 2
          \end{pmatrix}
        \)
        \label{fig:K5_15243}}
    \end{subfigure}\hfill
    \begin{subfigure}[t]{.3\linewidth}
        \centering
        \includegraphics[width=\linewidth]{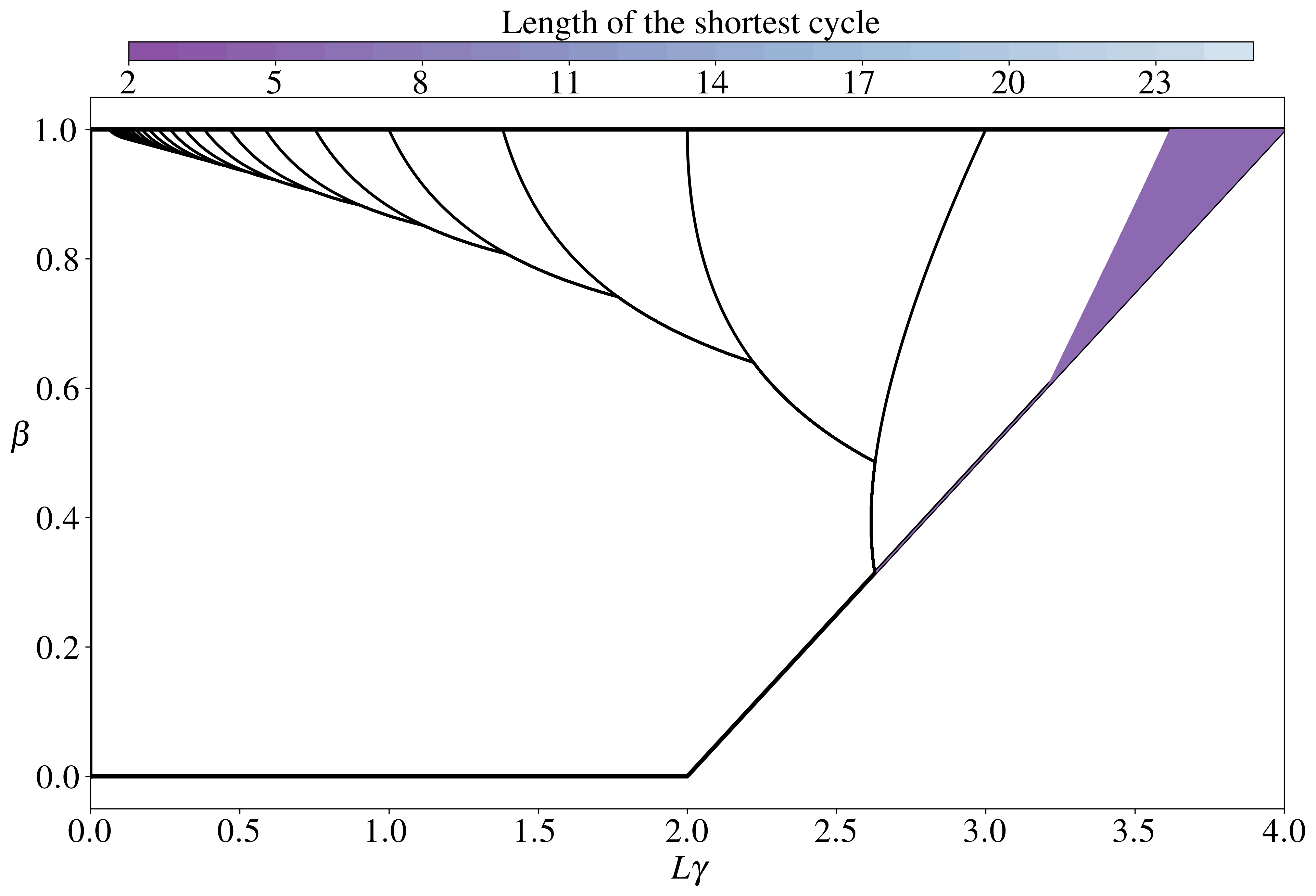}
        \caption{%
        \(\displaystyle 
          \sigma \;=\; 
          \begin{pmatrix}
            0 & 1 & 2 & 3 & 4 \\
            0 & 2 & 3 & 1 & 4
          \end{pmatrix}
        \)
        \label{fig:K5_13425}}
    \end{subfigure}\hfill
    
    \begin{subfigure}[t]{.3\linewidth}
        \centering
        \includegraphics[width=\linewidth]{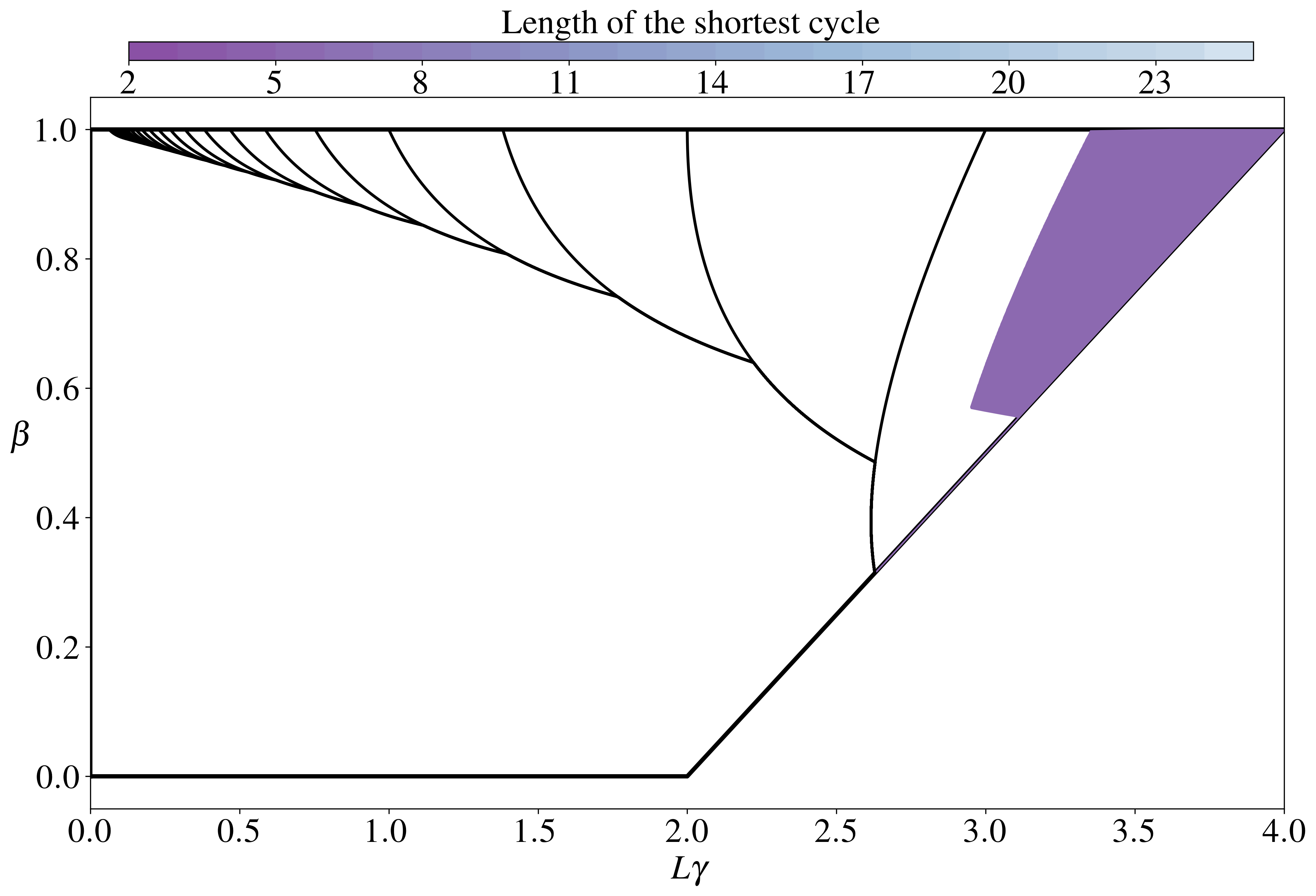}
        \caption{%
        \(\displaystyle 
          \sigma \;=\; 
          \begin{pmatrix}
            0 & 1 & 2 & 3 & 4 \\
            0 & 1 & 2 & 4 & 3
          \end{pmatrix}
        \)
        \label{fig:K5_12354}}
    \end{subfigure}\hfill
    \begin{subfigure}[t]{.3\linewidth}
        \centering
        \includegraphics[width=\linewidth]{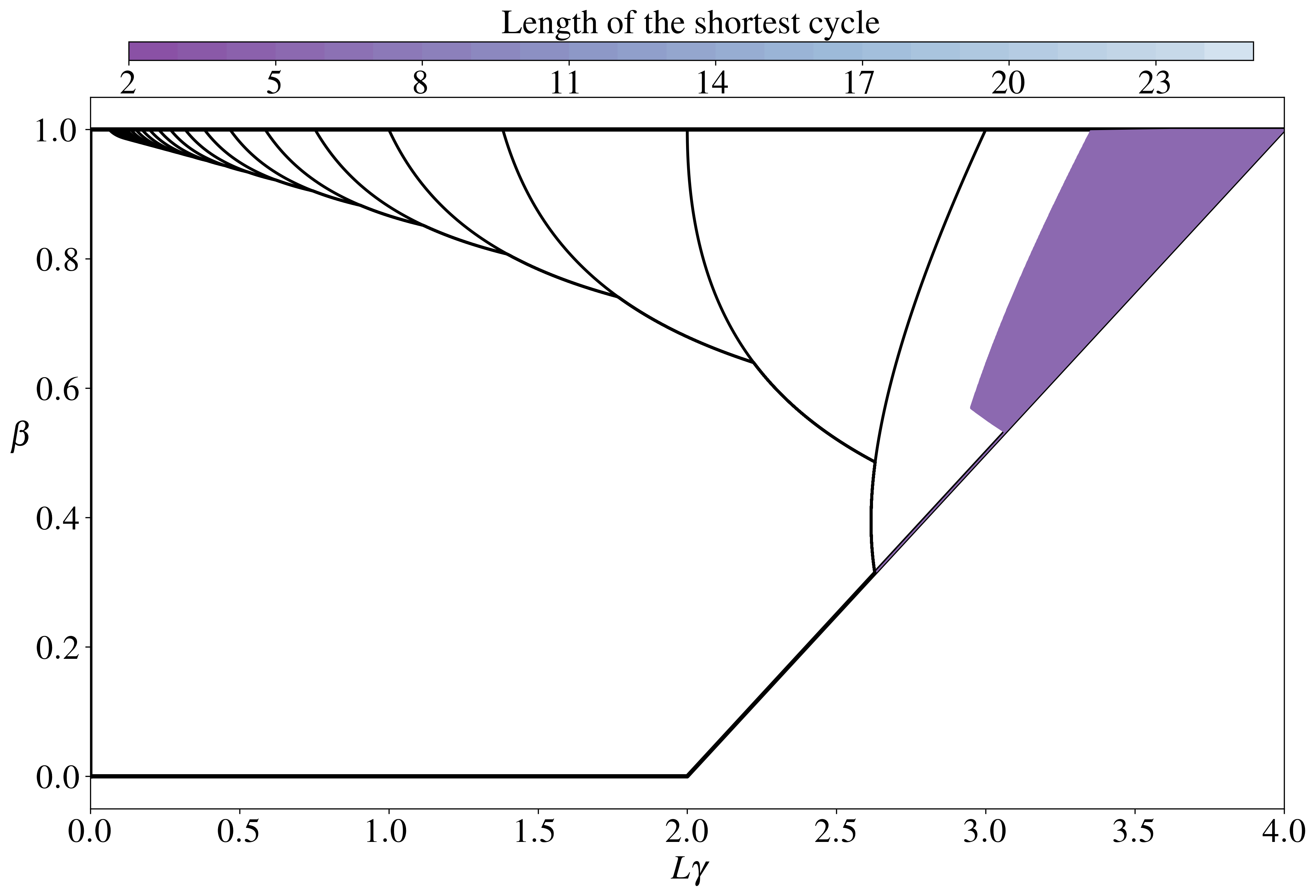}
        \caption{%
        \(\displaystyle 
          \sigma \;=\; 
          \begin{pmatrix}
            0 & 1 & 2 & 3 & 4 \\
            0 & 2 & 3 & 4 & 1
          \end{pmatrix}
        \)
        \label{fig:K5_13452}}
    \end{subfigure}\hfill
    \begin{subfigure}[t]{.3\linewidth}
        \centering
        \includegraphics[width=\linewidth]{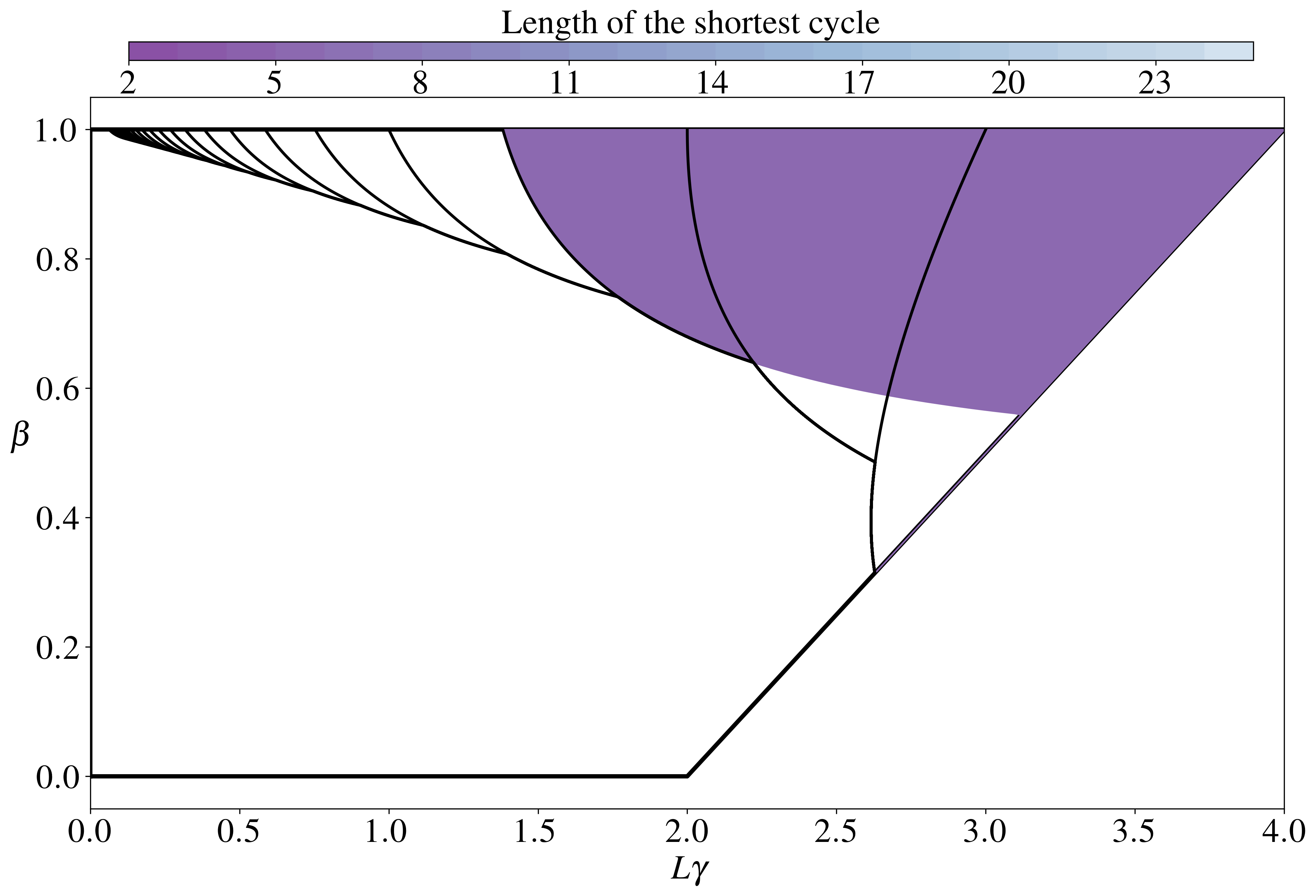}
        \caption{%
        \(\displaystyle 
          \sigma \;=\; 
          \begin{pmatrix}
            0 & 1 & 2 & 3 & 4 \\
            0 & 1 & 3 & 4 & 2
          \end{pmatrix}
        \)
        \label{fig:K5_12453}}
    \end{subfigure}
    
    \caption{For $K=5$, $(\mu,L)$ fixed, the set of points $\gb$ such that~\eqref{eq:lpsigma} admits a feasible point, depending on $\sigma$.
    For the 6 other permutations, the set is empty (as in e.g.~\Cref{fig:K4_1234}). \label{fig:cycles_dep_on_sig_K=5}}
\end{figure}

\paragraph{Identifying the correct permutation to obtain numerical results for $ K> 6$.} Having $K!$ possible permutations for a length $K$-cycle makes the search for long cycles untractable with the approach above. For example,  we could not have rendered with it the~\Cref{fig:cycle_d_2_1} with cycles of length 25. 
\begin{wrapfigure}[5]{r}{0.45\textwidth}
    \centering
    \includegraphics[width=.8\linewidth]{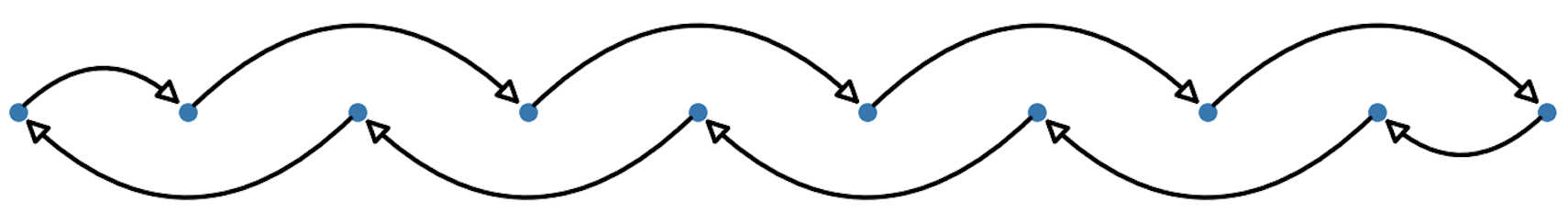}
    \caption{Typical permutation. \label{fig:permutation}}
\end{wrapfigure}
To go beyond $K=6$, we study the permutation for which cycles are obtained and identify a pattern, which enables us to conjecture that \textit{a single permutation} always enables to detect all cycles. More precisely, we conjecture that the one right permutation is such that $$(x_0, x_1, x_2, \cdots, x_{K-3}, x_{K-2}, x_{K-1}) = (x^{(0)}, x^{(1)}, x^{(3)}, x^{(5)}, \cdots, x^{(4)}, x^{(2)}),$$ as illustrated in \Cref{fig:permutation}. This conjecture comes from \Cref{fig:cycles_dep_on_sig_K=4,fig:cycles_dep_on_sig_K=5}, for $K=4$ and $K=5$ respectively.

\paragraph{Comments on \Cref{fig:cycles_dep_on_sig_K=4,fig:cycles_dep_on_sig_K=5}.}
To denote permutations, we write in two lines, the second line giving the images of the first line, i.e.~{\footnotesize $\begin{pmatrix}
0 & 1  & \dots & K-1 \\
\sigma(0) & \sigma(1) & \dots & \sigma(K-1)
\end{pmatrix}$}, e.g.~{\footnotesize$\sigma \;=\; 
\begin{pmatrix}
0 & 1 & 2 & 3 & \dots & K-2 & K-1 \\
0 & 1 & 3 & 5 & \dots & 4 & 2
\end{pmatrix}$} for the permutation of interest. Note that we limit the number of cycles to represent:
\begin{enumerate}
\item We can fix the first point of $X$ to be the smallest one (i.e., 0 is a fixed point of the permutation that orders points). We thus only have to consider $K-1!$ permutations (up to composition with the useless cyclic permutation $(0 \ 1 \, \cdots \, K-1)$).
\item By a symmetry argument, if there exists a feasible point for a permutation {\footnotesize $\begin{pmatrix}
0 & 1  & \dots & K-1 \\
\sigma(0) & \sigma(1)   & \dots & \sigma(K-1)
\end{pmatrix}$}, there exists a feasible point for the permutation {\footnotesize $\begin{pmatrix}
0 & 1 &\dots & K-1 \\
K-1-\sigma(0) & K-1-\sigma(1) & \dots & K-1-\sigma(K-1)
\end{pmatrix}$}, which reduces the number of permutation by a factor 2\footnote{(Almost! this is in fact only true for $K$ odd, for an even cycle length $K$, a couple of extra cases need to be considered,  when some permutations remain unchanged by this transformation, up to composition with the cyclic permutation used in the previous point, but this goes beyond this discussion).}.
\end{enumerate}

\paragraph{Conclusion on \Cref{fig:cycles_dep_on_sig_K=4,fig:cycles_dep_on_sig_K=5}.} 
We observe that for any tuning for which we could find a cycle, we could find one built with one of the permutations displayed in~\Cref{fig:permutation}. 
In other words, it appears that constraining the search to \textit{that single} permutation scheme enables obtaining the same purple area from~\Cref{fig:dim1}. We thus used this observation to draw~\Cref{fig:cycle_d_2_1} up to $K=25$.
Ultimately, we therefore have two connected conjectures:
\begin{enumerate}
    \item[(i)] Can we prove that the purple region representing all the cycles, is achieved in dimension 1, therefore proving that $\HB$ cannot accelerate in dimension 1 either?
    \item[(ii)] Can we prove it using only the cycle ordering given in \Cref{fig:permutation}?
\end{enumerate}

\paragraph{Link to \citet{wang2022provable}.}
While~\citet{wang2022provable} seem to show that $\HB$ accelerate in dimension 1, an incorrect inequality is used in Appendix E between the equations 84 and 85, as authors multiply an inequality with $1 - L\eta - 2\theta\eta$ as if the latter was non-negative. In other words, they use $1 - L\eta - 2\theta\eta \geq 0$, which cannot be true in their setting (equivalent to $0 < 2L \leq \mu \leq L$).

\section{Open problem 2: Lyapunov functions, cycles, and third type of behavior}
\vspace{-0.5em}
While green regions of~\Cref{fig:2} correspond to tunings on which $\HB$ admits a Lyapunov function (and hence converges), and purple regions of~\Cref{fig:2} correspond to tunings on which $\HB$ admits cycles (and hence does not converge in general), it seems that some tunings do not fall in any of the two latter regions. We can therefore ask whether $\HB$ converges (without any Lyapunov function) or diverges (without any cycling behavior) on those tunings. First, we provide some details on the form of the Lyapunov functions and cycles we looked for. Then we discuss the possible types of behavior that may occur in the remaining areas. Finally, we discuss a mindblowing observation that lacks an explanation.

\paragraph{More details on the Lyapunov functions we looked for.}

The green (including both shades of green) region of~\Cref{fig:2} has been drawn numerically following the strategy described in~\citep{taylor2018lyapunov}. We looked for a Lyapunov sequence $V_t$ (implicitly depending on $f$, $x_{t-1}$ and $x_t$) under the form $$V_t = \ell(f(x_t) - f_\star, f(x_{t-1})-f_\star) + q(x_{t-1} - x_\star, \nabla f(x_{t-1}), x_t-x_\star, \nabla f(x_t))$$ with $\ell$ a linear function and $q$ a homogeneously quadratic one (i.e.~linear combination of the inner products of pairs of inputs).

For a given tuning $\gamma, \beta$ of $\HB$, and a given rate $\rho$, we \emph{automatically} determine the existence of a Lyapunov sequence of the form described above and verifying
$$f(x_{t+1}) - f_\star \leq V_{t+1} \leq \rho V_t$$
for any sequence of iterates generated by $\HB_{\gb}$ applied on any function of $\Fml$.
The above-mentioned green region corresponds to the case where we found such a Lyapunov sequence working for $\rho=1$. Note that, for each tuning in the interior of the green region, there exists a Lyapunov sequence working for some $\rho<1$, constituting a proof of convergence of $\HB_{\gb}$ on $\Fml$ for any $\gb$ picked in this region.

To disprove convergence elsewhere, we decided to exhibit cycle behaviors of the $\HB$ algorithm.

\paragraph{How did we find cycles?}

\begin{wrapfigure}{r}{0.35\textwidth}
    \centering
    \includegraphics[width=\linewidth]{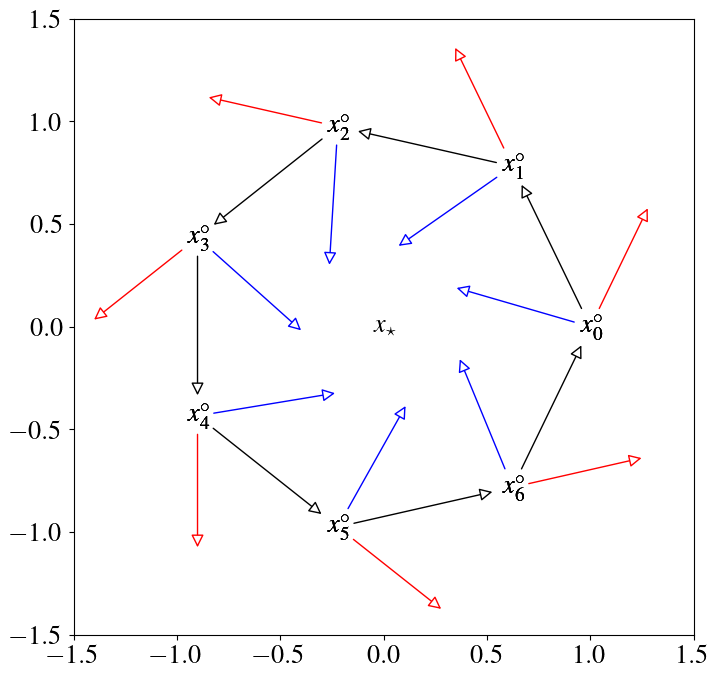}
    \caption{Cycle example. $\rightarrow$: algorithm updates. \textcolor{red}{$\rightarrow$}: momentum part. \textcolor{blue}{$\rightarrow$}: appropriate negative gradient.}
\end{wrapfigure}

Let $\gamma, \beta$ and $K\geq 3$ fixed.
Let us first define what we call a cycle of $\HB_{\gb}$ on $\Fml$.
A cycle of length $K$ is a $K$-uple $(x_0, \cdots, x_{K-1})$ such that there exists a function $f\in\Fml$ with $\HB_{\gb}$ generated sequence from $f$ being exactly $(x_{t \ \mathrm{mod} \ K})_{t \geq 0}$.

Prior to~\citep{goujaud2023provable}, we described a way to automate the search for such cycles using the PEP formalism in~\citep{goujaud2023counter}. It consists of simply running the algorithm and minimizing the distance between $(x_0, x_1)$ and $(x_K, x_{K+1})$ over the problem instances $f\in\Fml$. We show in~\citep{goujaud2023counter} that there is a cycle if and only if the infimum value of this optimization problem is 0. The advantage of this method is that we can use the PEP formalism described in~\citep{drori2014performance,taylor2017smooth} to formulate this as an SDP that can be solved efficiently. Moreover, there exist software (see~\citep{taylor2017performance,goujaud2022pepit}) to also automate the reformulation step to an SDP. Hence, determining the cycle existence of $\HB_{\gb}$ on $\Fml$ can be done in just a few of lines of code.

The issue with this first method is to determine whether the (numerically obtained) infimum value of this problem should be interpreted as a 0 with some machine precision error or as a positive value. This may lead to some uncertainty. To overcome this, we proposed another approach in~\citep{goujaud2023provable}.

\paragraph{What other type of behavior can we expect?}
Many other types of behavior could be possible. For instance, in~\citep[Future works section]{goujaud2023counter}, we provided an example of a function on which the vanilla gradient descent algorithm admits a chaotic behavior.

\paragraph{An interesting, unexplained observation.}

\begin{figure}[h!]
  \centering 
    \begin{subfigure}{.45\linewidth}
        \centering
        \includegraphics[width=\linewidth]{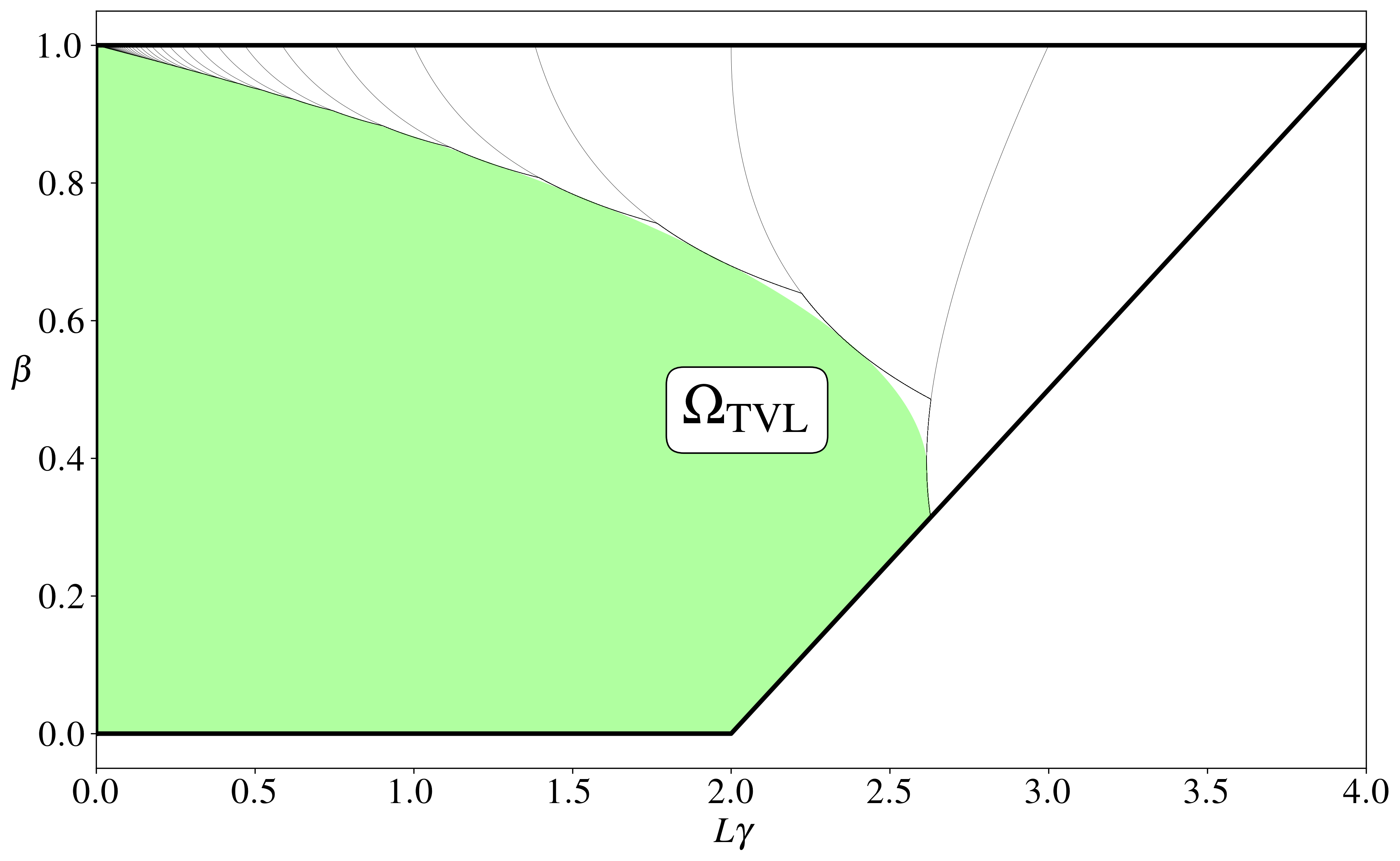}
        \caption{Cycle lengths multiples of 1 \label{fig:taylor_vs_rou_cycles}}
    \end{subfigure}
    \begin{subfigure}{.45\linewidth}
        \centering
        \includegraphics[width=\linewidth]{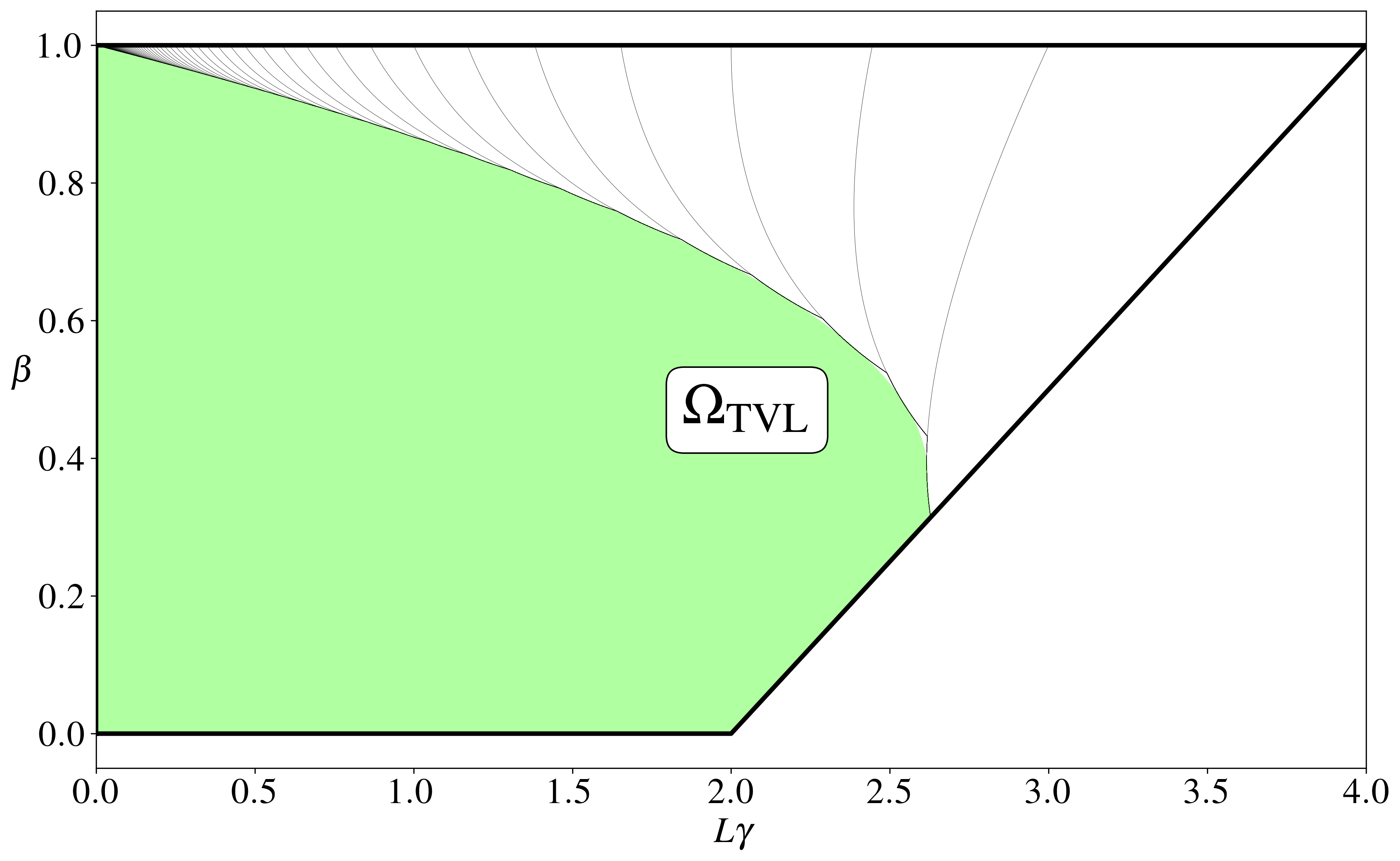}
        \caption{Cycle lengths multiples of 1/2}
    \end{subfigure}
    \begin{subfigure}{.45\linewidth}
        \centering
        \includegraphics[width=\linewidth]{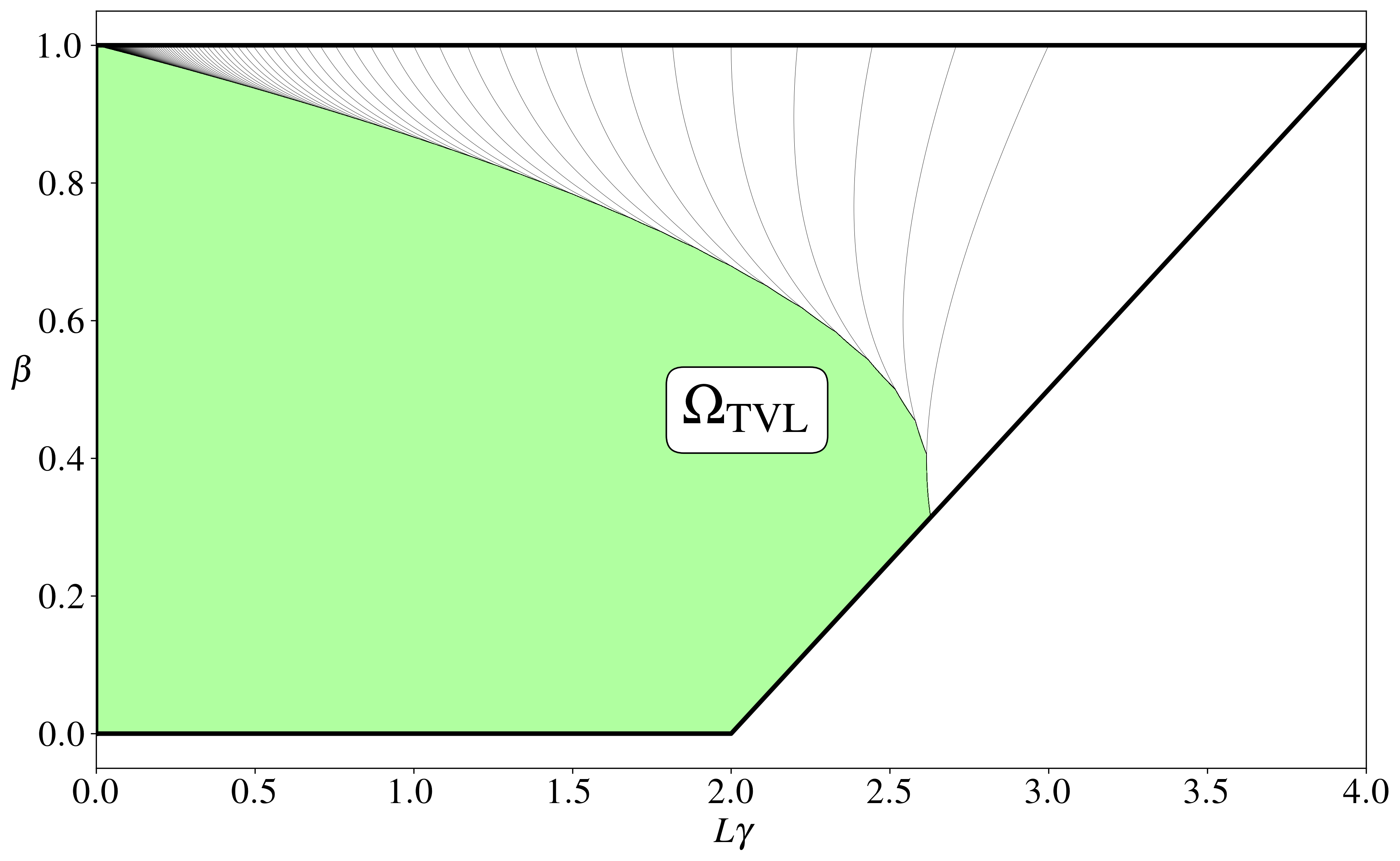}
        \caption{Cycle lengths multiples of 1/4}
    \end{subfigure}
    \begin{subfigure}{.45\linewidth}
        \centering
        \includegraphics[width=\linewidth]{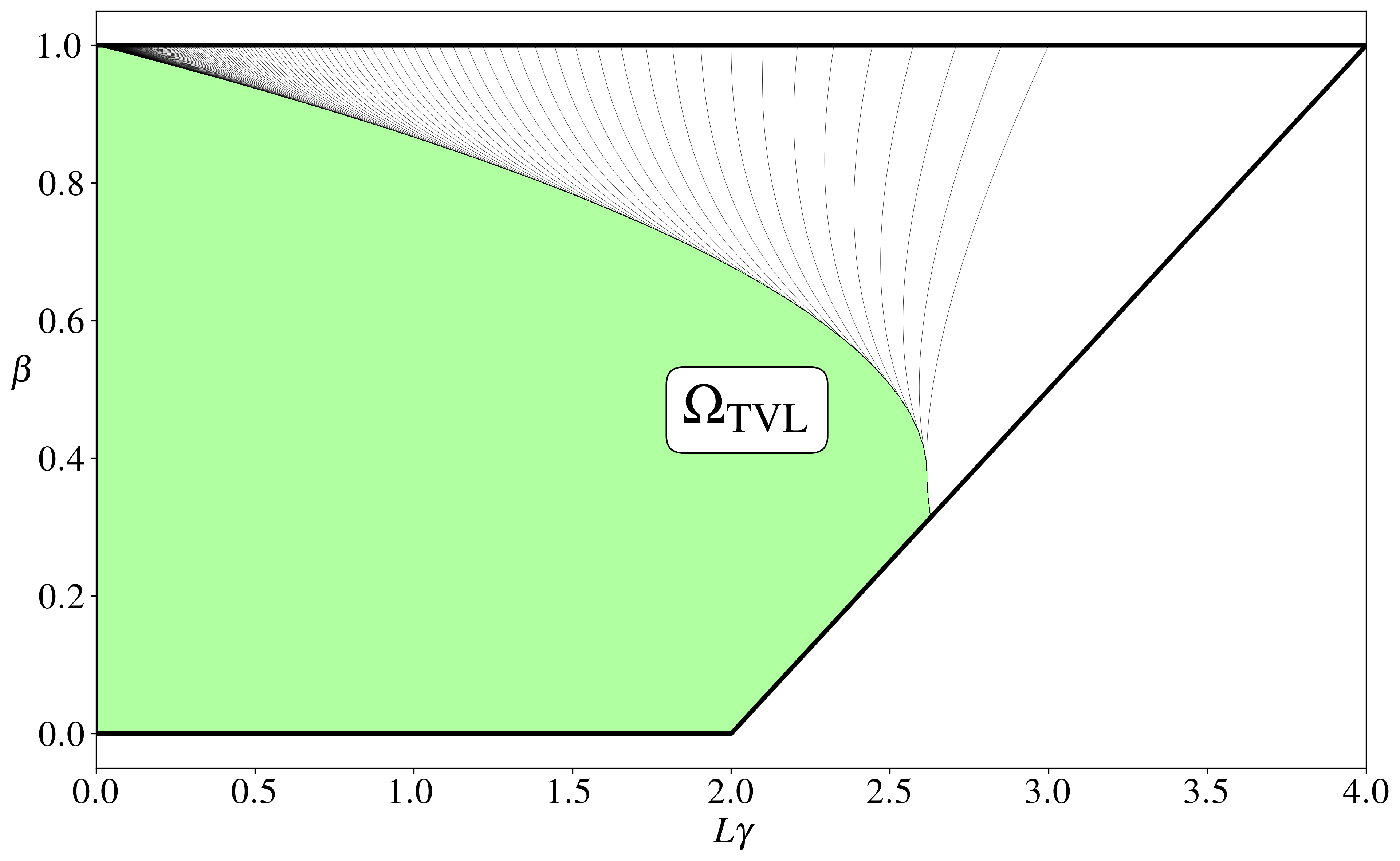}
        \caption{Cycle lengths multiples of 1/8}
    \end{subfigure}
    \caption{Area of Lyapunov in green. Cycles' borders are in black. Cycle lengths multiples of 1 are represented in the top left image, while cycle lengths multiple of 1/2, 1/4, and 1/8 are represented in three other ones. \label{fig:taylor_vs_rou_intermediate_cycles}}
\end{figure}

In~\Cref{fig:taylor_vs_rou_cycles}, we observe that the border of the Lyapunov region is smooth, while this is not the case for the cycling region. We see that this is due to the discrete union of smooth regions, each corresponding to a length of cycles.
Using the closed form formula of the border of each cycling region of a given length $K$ (that we obtained in~\citep[Theorem 3.5]{goujaud2023provable}), we displayed in~\Cref{fig:taylor_vs_rou_intermediate_cycles} some intermediate regions plugging fractional values of $K$ in the latest formula.

What we observe is that the union of a continuum of those regions seems to perfectly recover the blank regions without overlapping with the green ones, perfectly completing the figure.
We then naturally wonder: Since we did not define the notion of cycles of fractional length, what do those intermediate regions actually represent? What does it say about those blank regions?

\subsection{Acknowledgments}
The work of B. Goujaud and A. Dieuleveut is partly supported by ANR-19-CHIA-0002-01/chaire SCAI, and Hi!Paris FLAG project, PEPR Redeem.
A.~Taylor is supported by the European Union (ERC grant CASPER 101162889). Views and opinions expressed are however those of the author(s) only and do not necessarily reflect those of the European Union or the European Research Council Executive Agency (ERCEA). Neither the European Union nor the granting authority can be held responsible for them. The work of A. Taylor is also partly supported by the Agence Nationale de la Recherche as part of the ``France 2030'' program, reference ANR-23-IACL-0008 (PR[AI]RIE-PSAI).

\bibliographystyle{plainnat}
\bibliography{references}

\end{document}

%% file: header.tex
\usepackage{times}
\usepackage[utf8]{inputenc}
\usepackage[T1]{fontenc}
\usepackage{url}
\usepackage{booktabs}
\usepackage{multirow}
\usepackage{amsfonts}
\usepackage{yhmath}
\usepackage{nicefrac}
\usepackage{microtype}
\usepackage{faktor}
\usepackage{hhline}

\usepackage{thm-restate}

\usepackage{tikz}
\usepackage{graphicx}
\usepackage{grffile}
\usepackage{mathtools, stmaryrd}
\usepackage{tocloft}
\usepackage[hidelinks]{hyperref}
\usepackage{xcolor}
\hypersetup{
    colorlinks,
    linkcolor={red!50!black},
    citecolor={blue!50!black},
    urlcolor={blue!80!black}
}
\definecolor{ao(english)}{rgb}{0.0, 0.5, 0.0}

\usepackage[algo2e, algoruled, boxed, vlined]{algorithm2e}
\usepackage{algorithm}
\usepackage{algpseudocode}

\SetKwInput{KwInput}{Input}
\SetKwInput{KwOutput}{Output}

\usepackage[noabbrev]{cleveref}

\usepackage{caption}
\usepackage{subcaption}

\usepackage{amssymb}
\usepackage{amsmath}
\usepackage{natbib}
\usepackage{multirow}

\usepackage{pifont}

\usepackage{blindtext}
\usepackage{cancel}

\usepackage{enumitem}
\usepackage{float}
\usepackage{wrapfig}

\usepackage{tikz}

\usepackage{booktabs}
\usepackage{mdframed}
\newmdenv[topline=false,rightline=false]{leftbot}

\newtheorem{?}[Th]{Problem}

\crefname{Def}{Definition}{Definitions}
\crefname{Rem}{Remark}{Remarks}
\crefname{Th}{Theorem}{Theorems}

\newenvironment{proof}
{\ \\
\noindent \textit{Proof.}
}
{
$\hfill\blacksquare$}

\def\GD{\operatorname{GD}}
\def\HB{\operatorname{HB}}

\renewcommand{\leq}{\leqslant}
\renewcommand{\geq}{\geqslant}
\renewcommand{\le}{\leqslant}

\renewcommand{\epsilon}{\varepsilon}

\newcommand{\R}{\mathbb{R}}

\newcommand{\Qml}{\mathcal{Q}_{\mu, L}}
\newcommand{\Fml}{\mathcal{F}_{\mu, L}}

\newcommand{\cycle}{\mathrm{Cycle}}
\newcommand{\Ocyml}{\Omega_{\cycle}(\Fml)}

\newcommand{\Orouml}{\Omega_{\circ\text{-}\cycle}(\Fml)}

\newcommand{\gb}{\gamma,\beta}